\documentclass[11pt, leqno]{article}
\usepackage{}
\usepackage[total={6in, 9in}]{geometry}
 \usepackage{amsfonts}
\usepackage{amsthm}
\usepackage{amssymb}
\usepackage{amsmath,amsfonts}
\usepackage{mathrsfs}
\usepackage{cases}
\usepackage{latexsym,bm}
\usepackage{indentfirst}
\usepackage{color}
\usepackage{ifpdf}
\usepackage{graphicx}
\usepackage{psfrag}
\usepackage{dsfont}
\usepackage[ruled,vlined]{algorithm2e}
\usepackage{enumerate}
\usepackage[
pdfauthor={},
pdftitle={},
pdfstartview=XYZ,
bookmarks=true,
colorlinks=true,
linkcolor=blue,
urlcolor=blue,
citecolor=blue,
bookmarks=true,
linktocpage=true,
hyperindex=true
]{hyperref}

\usepackage{pgf,tikz,pgfplots}
\usepackage{tikz-3dplot}
\usepackage{pifont}
\usepackage{refcount}

\newtheorem{THM}{\textbf{Theorem}}[section]

\newtheorem{LEM}[THM]{\textbf{Lemma}}
\newtheorem{CLA}[THM]{\textbf{Claim}}
\newtheorem{CON}[THM]{\textbf{Conjecture}}

\newcommand{\pf}{\noindent\textbf{Proof}.\quad}

\newcommand{\ve}{\varepsilon }

\newcommand{\arxiv}[1]{\href{http://arxiv.org/abs/#1}{\texttt{arXiv:#1}}}

\begin{document}
\title{On the Multigraph Overfull Conjecture}
\author{Michael J. Plantholt and Songling Shan \\ 
	\medskip  
	Illinois State  University, Normal, IL 61790\\
	\medskip 
	{\tt mikep@ilstu.edu; sshan12@ilstu.edu}
}

\date{\today}
\maketitle

\emph{\textbf{Abstract}.}
 A subgraph $H$ of a multigraph $G$ is overfull if $ |E(H) | > \Delta(G) \lfloor |V(H)|/2 \rfloor.$
Analogous to  the  Overfull Conjecture proposed by Chetwynd and Hilton in 1986, Stiebitz et al. in 2012 formed  the  multigraph version  of the conjecture as follows:  Let $G$ be a multigraph with maximum multiplicity $r$ and maximum degree $\Delta>\frac{1}{3} r|V(G)|$. 
Then $G$ has chromatic index $\Delta(G)$ if and only if $G$ contains no overfull subgraph. 
In this paper, we prove the following three results toward the Multigraph Overfull Conjecture for sufficiently large and even $n$.

{\indent 
 (1) If $G$ is $k$-regular with $k\ge r(n/2+18)$, then $G$ has a 1-factorization. This    result  also settles a conjecture of the first author and 
 Tipnis from 2001 up to a constant error in the lower bound of $k$.

	(2) If $G$ contains an overfull subgraph and $\delta(G)\ge r(n/2+18)$, then $\chi'(G)=\lceil \chi'_f(G) \rceil$, where $\chi'_f(G)$ 
		is the fractional chromatic index of $G$.   
		
(3) If  the minimum degree of $G$ is  at least $(1+\ve)rn/2$ for any $0<\ve<1$ and $G$ contains no overfull subgraph, then $\chi'(G)=\Delta(G)$. 

}

The proof is based on the decomposition of multigraphs into simple graphs and we prove a slightly  weak version of a conjecture  due to the first author and Tipnis from  1991 on decomposing a multigraph into constrained  simple graphs. The result is of independent interests.   


\vspace{2mm}

\section{Introduction}

In this paper,  we use the term ``graphs''  for  multigraphs, which may contain multiple edges but contain no loop. 
A multigraph with no parallel  edge will be stressed as a simple graph. 
Let $G$ be a graph. 
Denote by $V(G)$ and  $E(G)$ the vertex set and edge set of $G$,
respectively, and by $e(G)$ the cardinality of $E(G)$. 
For $v\in V(G)$, $N_G(v)$ is the set of neighbors of $v$ 
in $G$, and 
$d_G(v)$, the degree of $v$
in $G$, is the number of edges of $G$ that are incident with $v$.
The notation $\delta(G)$ and $\Delta(G)$ denote the minimum degree 
and maximum degree of $G$, respectively.

For two integers $p,q$, let $[p,q]=\{ i\in \mathbb{Z} \,:\, p \le i \le q\}$.  
For an integer $k\ge 0$, an \emph{edge $k$-coloring} of a multigraph $G$ is a mapping $\varphi$ from $E(G)$ to the set of integers
$[1,k]$, called \emph{colors}, such that  no two adjacent edges receive the same color with respect to $\varphi$.  
The \emph{chromatic index} of $G$, denoted $\chi'(G)$, is defined to be the smallest integer $k$ so that $G$ has an edge $k$-coloring.   A color class of $\varphi$  is a set of edges of $G$ colored by
the same color under $\varphi$, which is a matching.

In the 1960's, Gupta~\cite{Gupta-67}  and, independently, Vizing~\cite{Vizing-2-classes}  proved
 that for all graphs $G$,  $\Delta(G) \le \chi'(G) \le \Delta(G)+\mu(G)$, 
 where $\mu(G)$, called the \emph{maximum multiplicity} of $G$, is the largest number of 
 edges joining two vertices of $G$.   Thus if $G$ is simple,  $\chi'(G)$ is always equal to either $\Delta(G)$ or $\Delta(G)+1$.  For multigraphs, the range of possible values for $\chi'(G)$ is wider.  Nevertheless, we show in this paper that for multigraphs of large order $n$ and minimum degree a bit more than $\mu(G) n/2$, we can generally determine the chromatic index exactly. This method is based on the concept of overfullness. 
 We say a  graph $G$  is \emph{overfull}
 if $|E(G)|>\Delta(G) \lfloor |V(G)|/2\rfloor$. 
    A  subgraph $H$ of $G$ is a \emph{$\Delta(G)$-overfull} subgraph
    if  $\Delta(H)=\Delta(G)$ and $H$ is overfull. 
    A subgraph $H$ of odd order of $G$ is \emph{$\Delta(G)$-full} 
    if $\Delta(H)=\Delta(G)$ and $|E(H)|=\Delta(H) \lfloor |V(H)|/2\rfloor$. 
    
  It is clear  that if $G$ contains a $\Delta(G)$-overfull subgraph, then 
 $\chi'(G) \ge \Delta(G)+1$.
  Conversely, for simple graphs, 
  Chetwynd and  Hilton~\cite{MR848854,MR975994},  in 1986, proposed the following 
conjecture. 
\begin{CON}[Overfull Conjecture]\label{overfull-con}
	Let $G$ be a simple graph  with $\Delta(G)>\frac{1}{3}|V(G)|$. Then $\chi'(G)=\Delta(G)$  if and only if $G$ contains no $\Delta(G)$-overfull subgraph.  
\end{CON}

The degree condition  $\Delta(G)>\frac{1}{3}|V(G)|$ in the conjecture above is best possible as seen by 
the simple graph $P^*$ obtained from the Petersen graph by deleting one vertex. 
 Applying Edmonds' matching polytope theorem, Seymour~\cite{seymour79}  showed  that whether a graph  $G$ contains an overfull subgraph of maximum degree $\Delta(G)$ can be determined in polynomial time. Thus if the Overfull Conjecture is true, then the NP-complete problem of 
determining the chromatic index~\cite{Holyer} becomes  polynomial-time solvable 
for simple graphs $G$ with $\Delta(G)>\frac{1}{3}|V(G)|$.
There have been some fairly strong results supporting the Overfull Conjecture in the case when $G$ is regular or $G$ has 
large minimum degree, for example, see~\cite{MR975994, MR3545109,MR1439301,MR4394718, 2105.05286, 2205.08564}. 
Our goal in this paper is to study the multigraph version~\cite[p. 259]{StiebSTF-Book} of the Overfull Conjecture first formed by Stiebitz et al. in 2012, which can be stated as follows.

\begin{CON}[Multigraph Overfull Conjecture]\label{overfull-con2}
	Let $G$ be a graph  satisfies $\Delta(G)>\frac{1}{3} \mu(G)|V(G)|$. Then $\chi'(G)=\Delta(G)$  if and only if $G$ contains no $\Delta(G)$-overfull subgraph.  
\end{CON}

Again, the degree condition  $\Delta(G)>\frac{1}{3}\mu(G)|V(G)|$ in Conjecture~\ref{overfull-con2} is best possible.  
To see this, let $r\ge 2$
be an  integer, and 
let $Q$ be obtained from the Petersen graph by duplicating each of 
its edge $r-1$ times. Then $Q$ is $3r$-regular with maximum multiplicity $r$. 
Let $Q^*$ be obtained from $Q$ by deleting a vertex.  It is easy to see that $Q^*$
contains no $3r$-overfull subgraph.  By a result of the first author and Tipnis~\cite[Theorem 2]{MR1483445}, 
we have  $\chi'(Q^*)=3r+1$ if $r$ is odd. Thus 
$Q^*$ is a sharpness example for the condition $\Delta(G)>\frac{1}{3} \mu(G)|V(G)|$.  (When $r$ is even, $\chi'(Q)=\chi'(Q^*)=3r$.)

First consider the case of regular graphs.
It is easy to verify that when $G$ is regular with even order, $G$ has no $\Delta(G)$-overfull subgraphs if its vertex degrees are at least about $\frac{1}{2}|V(G)|$.  Thus the well-known 1-Factorization Conjecture,  first stated in~\cite{MR772711} but may go back to Dirac in the early 1950s,  is a special case of the Overfull Conjecture.

\begin{CON}[1-Factorization Conjecture]\label{con:1-factorization}
	Let  $G$ be a simple graph  of even order $n$. If $G$ is $k$-regular for some  $k\ge 2\lceil n/4\rceil-1$,  then  $G$ is 1-factorable; equivalently, $ \chi'(G) = \Delta(G)$.
\end{CON}

In 2016,  Csaba, K\"uhn, Lo, Osthus and Treglown~\cite{MR3545109} verified Conjecture~\ref{con:1-factorization} for sufficiently large $n$. 

\begin{THM}\label{thm:1-factorization-proof}
	There exists an $n_0\in \mathbb{N}$ such that the following holds. Let $n, k\in \mathbb{N}$
	be such that $n\ge n_0$ is even and  $k\ge 2\lceil n/4\rceil-1$.  
	Then every  $k$-regular  simple graph on  $n$  vertices has a 1-factorization.
\end{THM}

A natural extension of the 1-factorization Conjecture to multigraphs is obtained by restricting the edge multiplicity.  The conjecture 
was  proposed by the first author and 
Tipnis in  2001~\cite{MR1877660}; see also in the ``Graph Edge Coloring'' book ~\cite[p. 260]{StiebSTF-Book}.

\begin{CON}[Multigraph 1-Factorization Conjecture]\label{con:1-factorization-multi}
Let  $G$ be a graph  of even order $n$ and maximum multiplicity $r$. If $G$ is $k$-regular for some  $k\ge  r(2\lceil n/4\rceil-1)$,  then  $G$ is 1-factorable; equivalently, $ \chi'(G) = \Delta(G)$.
\end{CON}

Conjecture~\ref{con:1-factorization-multi}  was proved when  $k$ is large by the first author and Tipnis~\cite{MR1149003, MR1877660}. 
Vaughan~\cite{MR2993074} proved Conjecture~\ref{con:1-factorization-multi} 
for large graphs asymptotically.  
In this paper, we prove a slightly weaker form of Conjecture~\ref{con:1-factorization-multi}.  

\begin{THM}\label{thm:1-factorization}
There exists an $n_0\in \mathbb{N}$ such that the following holds. Let $n, k, r\in \mathbb{N}$
be such that $n\ge n_0$ is even and  $k\ge r(n/2+18)$. Then every $k$-regular graph $G$ on
$n$ vertices  with maximum multiplicity at most $r$ has a 1-factorization. 
\end{THM}

For a 
graph $H$ with an odd number of vertices $n\ge 3$,
we have $\chi'(H) \ge  \frac{|E(H)|}{(n-1)/2}$, since each color 
class is a matching and each matching can contain at most $(n-1)/2$ edges.  
For any graph $G$ on at least three vertices, define the \emph{density} of $G$ as 
$$
\omega(G)=\max\left\{\frac{|E(H)|}{(n-1)/2}: \text{$H\subseteq G$, $|V(H)|=n$, $n\ge 3$ and is odd}\right\}. 
$$
It is clear that $\chi'(G) \ge \omega(G)$. Combining with the 
lower bound $\Delta(G)$ on $\chi'(G)$, 
$\chi_f'(G):=\max\{\Delta(G), \omega(G)\}$ 
is called the \emph{fractional chromatic index} of $G$. As the chromatic index is always an 
integer, we have a general lower bound for the chromatic index given by 
$\chi'(G) \ge \lceil \chi_f'(G) \rceil$. 
We  show that when $\Delta(G)$-overfullness or  $\Delta(G)$-fullness is present, then the chromatic index  of $G$ equals the integer round-up of the fractional chromatic index.

\begin{THM}\label{thm:overfull-present}
	There exists an $n_0\in \mathbb{N}$ such that the following holds. Let $n,  r\in \mathbb{N}$
	be such that $n\ge n_0$ is even.  
	If $G$ is a graph of order $n$,  maximum multiplicity $r$, $\delta(G) \ge r(n/2+18)$,   and $G$ contains a $\Delta(G)$-full 
	or $\Delta(G)$-overfull subgraph, then $\chi'(G)=\lceil \chi_f'(G)\rceil$.  
	As a consequence, $\chi'(G)=\Delta(G)$ if  $G$ contains a $\Delta(G)$-full subgraph. 
\end{THM}

Restricting the Overfull Conjecture on simple graphs of even order and large minimum degree, Theorem~\ref{thm:plantholt-shan} was proved by the authors. In this paper, 
we also prove its analogy for multigraphs. 

\begin{THM}[{\cite{2105.05286}}]\label{thm:plantholt-shan}
	For all $0<\ve <1$, there exists an $n_0\in \mathbb{N}$ such that the following holds. Let $n \in \mathbb{N}$
	be such that $n\ge n_0$ is even. 
	If $G$ is a simple graph on $n$ vertices with $\delta(G) \ge (1+\ve)n/2$,   then $\chi'(G)=\Delta(G)$ if and only if $G$  contains no $\Delta(G)$-overfull subgraph. 
\end{THM}

\begin{THM}\label{thm:1}
	For all $0<\ve <1$, there exists an $n_0\in \mathbb{N}$ such that the following holds. Let $n, r \in \mathbb{N}$
	be such that $n\ge n_0$ is even.  
	If $G$ is a graph on $n$ vertices with maximum multiplicity $r$ and $\delta(G) \ge r(1+\ve)n/2$,   then $\chi'(G)=\Delta(G)$ if and only if $G$  contains no $\Delta(G)$-overfull subgraph. 
\end{THM}

Generally speaking, the proofs of Theorems~\ref{thm:1-factorization},  \ref{thm:overfull-present} and ~\ref{thm:1} are based on the decompositions of multigraphs into simple graphs  with constrained minimum degree conditions and 
then an application of the existing  results  Theorem  \ref{thm:1-factorization-proof}  or Theorem~\ref{thm:plantholt-shan}  on  each of the simple graphs from the decomposition.    The hardness in the decomposition is to  avoid the overfullness in each of the simple graph from the decomposition and have the simple graph satisfying degree constrains. 
In doing so, we develop a decomposition result for regular  multigraphs  of large  degree and  odd maximum multiplicity (see Lemma~\ref{lem:decompose2}),   which 
solves a slightly weak version of a conjecture due to
 the first author and Tipnis~\cite{MR1149003}  from 1991.

The remainder of this paper is organized as follows.
In the next section, we introduce some notation and   preliminary results.
In   Section 3, we prove Theorems~\ref{thm:1-factorization} and~\ref{thm:overfull-present}, and in Section 4,
we prove Theorem~\ref{thm:1}.

\section{Notation and preliminaries}

Let $G$ be a graph 
 and $A,
B\subseteq V(G)$ be two disjoint vertex sets. Then $E_G(A,B)$ is the set
of edges in $G$  with one end in $A$ and the other end in $B$, and  $e_G(A,B):=|E_G(A,B)|$.  We write $E_G(v,B)$ and $e_G(v,B)$
if $A=\{v\}$ is a singleton.  For an edge $e\in E_G(u,v)$, if $e_G(u,v)=1$, we call $e$ a \emph{singleton edge}. 
For 
$S\subseteq V(G)$,   
the subgraph of $G$ induced by  $S$ is  $G[S]$, and  $G-S:=G[V(G)\setminus S]$. 
If $F\subseteq E(G)$, then $G-F$ is obtained from $G$ by deleting all
the edges of $F$.  Denote by $V_\Delta$ the set of maximum  degree vertices of $G$.

A \emph{trail} is an alternating sequence of vertices and edges $v_0e_1v_1\ldots e_tv_t$ such that
$v_{i-1}$ and $v_i$ are the end vertices of $e_i$
for each $i\in [1,t]$, and the edges are all distinct (but there might be repetitions among the vertices).  A trail is \emph{closed} if $v_0= v_t$, and is \emph{open} otherwise. An {\it Euler tour} of $G$ 
is a closed trail in $G$ that contains all the edges of $G$.  
A graph is \emph{even} if all its vertex degrees are even. 
We will need the following classic result of Euler. 
\begin{THM}[Euler, 1736]\label{Euler}
	A graph $G$ has an Euler tour if and only if 
	$G$ has at most one nontrivial component and $G$ is even. 
\end{THM}

\begin{LEM}\label{cor:cycle-path-decomp}
If $G$ is a graph with in total $2\ell$ vertices of odd degree for some integer $\ell\ge 0$, 
then $G$ can be decomposed into  edge-disjoint cycles and paths, where there are exactly $\ell$
paths in the decomposition such that the set of the endvertices of the paths 
is the same as the set of odd degree vertices of $G$ and that the union of the paths is a forest.   
\end{LEM}

\pf  First  we iteratively delete the edges of cycles from $G$ until no cycles are left.  The remaining graph $G^*$ is a forest, and has the same set of odd degree vertices as does $G$.  From a component of $G^*$, remove a path between two endvertices; this reduces the number of odd degree vertices by two.  Then iterating this to get  $\ell$ paths, we get a decomposition with the desired properties.
\qed 
 
 We will also need the two classic results below on Hamilton cycles in simple graphs. 
 
 \begin{THM}[\cite{MR47308}]\label{thm:Dirac}
 	If  $G$ is a simple graph on $n\ge 3$ vertices with $\delta(G) \ge \frac{n}{2}$, then $G$ has a Hamilton cycle. 
 \end{THM}

 \begin{THM}[\cite{MR294155}] \label{lem:chvatal's-theorem}
 	Let $G$ be a simple graph on  $n\ge 3$ vertices. 
 	Suppose the degrees of  $G$ are $d_1, \ldots d_n$ with $d_1\le \ldots \le d_n$. 
 	If $d_i\ge i+1$ or $d_{n-i} \ge n-i$
 	for all $i<\frac{n}{2}$, then $G$ has a Hamilton cycle. 
 \end{THM}


 Let $f$ be a function from the vertices of a graph $G$ into the positive integers, $g$ be a map from possible edges into the  positive integers.  An  \emph{$fg$-coloring} of $G$ is a coloring of the edges so that each vertex $v$  has at most $f(v)$ incident edges assigned the same color, and  for each pair $u,v\in V(G)$, there are   at most  $g(uv)$ edges of the same color joining $u$ and $v$.  The results below provide an upper 
 bound on the number of colors needed for an $fg$-coloring. 
 
  \begin{LEM}[{\cite[Theorem 8]{MR1075068}}]\label{lem:fg-b}
Every bipartite graph $G$ has an $fg$-coloring in  $$\max \{\max_{v\in V(G)} \lceil d_G(v)/f(v) \rceil, \max_{u,v\in V(G)} \lceil e_G(u,v)/g(uv)\rceil\}$$ colors.
 \end{LEM}
 
 \begin{LEM}[{\cite[Corollary 15]{MR1075068}}]\label{lem:fg}
 If $f(v)$ and $g(uv)$ are all positive even integers for each vertex $v$ of 
$G$ and each pair $u,v$ of $V(G)$, then  $G$ has an $fg$-coloring in  $$\max \{\max_{v\in V(G)} \lceil d_G(v)/f(v) \rceil, \max_{u,v\in V(G)} \lceil e_G(u,v)/g(uv)\rceil\}$$ colors.
  \end{LEM}

  \begin{LEM}\label{lem:numberofD}
  	Let $G$ be an $n$-vertex  graph with maximum multiplicity $r$ such that all  vertices of degree less than $\Delta(G)$ are mutually adjacent in $G$ with $r$ edges between them. Then $|V_\Delta|> \frac{n}{2}$. 
  \end{LEM}
  
  \pf  By contradiction. Let $A=V_\Delta$, $|A|=k$ for some integer $k\ge 1$, and suppose that  $|V(G)\setminus A|=k+s$
  for some integer $s\ge 0$.  Let $H$ be obtained 
  from $G$ by deleting $s$ vertices from $V(G)\setminus A$, and $B=V(H)\setminus A$. Then 
  we have $|A|=|B|=k$.  
  As all the vertices of $B$ are mutually adjacent in $G$ with $r$ edges between them, 
  we have $\sum_{v\in B}d_H(v)=2r{k \choose 2}+e_H(A,B) \ge \sum_{v\in A}d_H(v)$. 
  Since $d_G(v)=d_H(v)+rs$ for each $v\in B$ and $d_G(v) \le d_H(v)+rs$ for each $v\in A$,
  we then get 
  $\sum_{v\in B}d_G(v)\ge \sum_{v\in A}d_G(v)$,  
  a contradiction. 
  \qed


\begin{LEM}\label{lemma:overfull-subgraph2}
	Let $G$ be a graph of even order $n$. There is at most one vertex $v$,
of minimum degree,  
	such that $G-v$ is $\Delta(G)$-overfull; if $G$ is not regular, there are at most two 
	vertices $v, w$, both of minimum degree,  such that $G-v$ and $G-w$ are $\Delta(G)$-full. 
\end{LEM}

\pf  Let   $v, w\in V(G)$. 
Suppose $G-v$ is $\Delta(G)$-full or $\Delta(G)$-overfull. 
 Note that  $e(G) \le( \Delta(G) (n-2)+d_G(v)+d_G(w))/2$. Thus 
\begin{equation}\label{eqn:overfull}
	e(G-v)=e(G)-d_G(v) \le \frac{\Delta(G) (n-2)-d_G(v)+d_G(w)}{2}.
\end{equation}
From the inequality above, we have $d_G(w)>d_G(v)$ if 
 $G-v$ is $\Delta(G)$-overfull, i.e., $2e(G-v)/(n-2)>\Delta(G)$. 
 Therefore, $v$ is the only vertex of minimum degree in $G$. 
 Similarly, suppose $G-v$ is $\Delta(G)$-full. Then we have $d_G(w) \ge d_G(v)$
 and so $v$ is a vertex of minimum degree in $G$.   If $G$ is not regular and 
 $d_G(w)=d_G(v)$, we must have $d_G(u)=\Delta(G)$
 for any $u\in V(G)\setminus\{v,w\}$ by~\eqref{eqn:overfull}. 
 Therefore,  if $G$ is not regular, there are at most two 
 vertices $v, w$, both of minimum degree,  such that $G-v$ and $G-w$ are $\Delta(G)$-full. 
 \qed 

In the case of regular graphs of high degree, it is easy to show there can be no overfull subgraphs.  The next lemma shows that  overfullness has a rigid structure when the minimum degree is high.

\begin{LEM}\label{lemma:overfull-subgraph}
	Let $G$ be a graph of even order $n$ and $\delta(G)> r\frac{n}{2}$, where $r=\mu(G)$.  
	Then every $\Delta(G)$-full or $\Delta(G)$-overfull subgraph is 
	obtained from $G$ by deleting a vertex of minimum degree. 
\end{LEM}

\pf Suppose to the contrary that $G[X]$ is  $\Delta(G)$-full or  $\Delta(G)$-overfull for some $X\subseteq V(G)$
with  $|X|$ odd and $3\le |X| \le n-3$.  Then $e(G[X])  \ge \Delta(G) \frac{|X|-1}{2} >r\frac{n}{2} \frac{|X|-1}{2}$. 
Since $e(G[X]) \le \frac{r|X|(|X|-1)}{2}$, we get $|X| >n/2$. 
 As $e(G[X]) \le  r\frac{|X|(|X|-1)}{2}$ and so $\omega(G[X]) \le r|X|$, 
 we get $\Delta(G) \le r|X|$. 
 As $ 3\le |V(G)\setminus X| =n-|X|\le \frac{n-1}{2}$ and $n$ is even, we get $n\ge 8$. 
 Now  since $(n-|X|)(\delta(G)-r(n-|X|-1))-r|X|$ is a quadratic concave down function in $|X|$
 and so its minimum is achieved at the two boundary values of $|X|$, we get 
 \begin{eqnarray*}
 	e_G(V(G)\setminus X, X)-r|X|  &\ge&  (n-|X|)(\delta(G)-r(n-|X|-1))-r|X|\\
 	 &\ge & \min\{3\delta(G)-6r-r(n-3), \frac{3r}{2}(n-1)/2-r(n+1)/2\}  \\
 	 & >& 0. 
 \end{eqnarray*}
Thus $e_G(V(G)\setminus X, X) > r|X|  \ge \Delta(G)$. 
This shows that $G[X]$ is  neither $\Delta(G)$-full nor $\Delta(G)$-overfull, a contradiction. 
The
lemma  then follows easily by Lemma~\ref{lemma:overfull-subgraph2} and
the fact that 
 $e(G-u)\le e(G-v)$ for any $u,v\in V(G)$ with $d_G(v)=\delta(v)$.  
 \qed

 \begin{LEM}[{\cite[Lemma 1]{MR1149003}}]\label{lem:spanning forest}
 If $G$ is a connected even graph, then $G$ has a spanning tree $T$
 such that $d_T(v) \le 1+\frac{1}{2}d_G(v)$ for each $v\in V(G)$. 
\end{LEM}
\section{Proof of Theorems~\ref{thm:1-factorization} and~\ref{thm:overfull-present}}

We need the following result by the first author and  Tipnis to 
decompose a graph into simple graphs. 

\begin{LEM}[{\cite[Theorem 1]{MR1149003}}]\label{lem:decomposition}
	Let  $G$  be a  $k$-regular graph with  even order $n$ and $k= rs \ge r(n/2 + 1)$, where $r=\mu(G)$ and $s$ is a positive integer.  If $r$ is even,  
	then $G$ can be decomposed into $ r/2$  edge-disjoint Hamilton cycles  and $r$  spanning simple graphs, each of which is  $(s-1)$-regular.
\end{LEM}

For  the graph $G$ described above,  when $\mu(G)$ is odd, the first author and Tipnis~\cite{MR1149003} 
conjectured that $G$ can be decomposed into $\mu(G)$ perfect matchings and   $\mu(G)$ simple graphs that each are $(s-1)$-regular. We here prove a weaker version of the conjecture, which will also be used to prove Theorem~\ref{thm:1-factorization}.  

\begin{LEM}\label{lem:decompose2}
Let  $G$  be a  $k$-regular graph with  even order $n\ge 96$ and  odd maximum multiplicity $r\ge 3$. Then the following statements 
hold. 
\begin{enumerate}[(a)]
	
	\item If $k= r s \ge r(n /2+ 16)$   for some even integer $s$,   then 
	$G$ can be decomposed into $r-2$ perfect matchings,    $(r+1)/2$ 
	Hamilton cycles, one spanning  $(s-1)$-regular simple graph, and $(r-1)$  spanning $(s-2)$-regular simple graphs.  
	\item If $k = rs \ge r(n/2+ 17)$ for some odd integer $s$,   then 
	$G$ can be decomposed into $2r-2$ perfect matchings,    $(r+1)/2$ 
	Hamilton cycles, one  spanning $(s-2)$-regular simple graph, and $(r-1)$  spanning $(s-3)$-regular simple graphs.  
\end{enumerate}
\end{LEM}

\pf  For (b),  the conclusion follows  from Statement (a) by  deleting  $r$  edge-disjoint perfect matchings from $G$. 
Thus we only need to prove Statement (a). 
Our strategy is to double each edge of $G$ to get a  $k^*$-regular graph $G^*$
with $k^* = 2r s \ge 2r(n/2 + 16)$ and then ``split off''
one Hamilton cycle $C$ of $G^*$ and one  $(s-1)$-regular 
subgraph $H_0$ of $G^*$ such that these two graphs are edge-disjoint and 
they both are subgraphs of $G$ and $\mu(G-E(C\cup H_0)) \le r-1$. 
Since  $G-E(C\cup H_0)$ is $((r-1)s-1)$-regular, then we can take off 
$s-2$ edge-disjoint perfect matchings from $G-E(C\cup H_0)$, and 
apply Lemma~\ref{lem:decomposition} on the remaining graph.

We double all the edges of $G$ to get $G^*$. Then $G^*$ is $ 2r s$-regular with maximum multiplicity $2r$.
We claim that $G^*$ has an orientation $D$
such that 
\begin{enumerate}[(1)]
	\item  $d_D^+(v)=d_D^-(v)$ for any $v\in V(D)$;  
	\item for any $u,v\in V(D)$, there are at most $r$ arcs of $D$ from $u$ to $v$; 
\end{enumerate}

To get such an orientation of $G^*$, for each edge $e$ of $G$ that  is not a singleton edge,
we pair up $e$ with its duplication  $e'$ and orient them in opposite  directions.
We let $R$ be the subgraph of $G$ induced by all its singleton edges.
Suppose $R$ has in total $2\ell$ vertices of odd degree for some integer $\ell \ge 0$. By Lemma~\ref{cor:cycle-path-decomp}, 
 $R$ can be decomposed into edge-disjoint cycles, say $C_1, \ldots, C_m$, together 
 with $\ell$ edge-disjoint 
paths $P_1, \ldots, P_\ell$, where the set of the endvertices of those paths  is the same as the set of odd degree vertices of $R$
and the union of these paths is a forest. 
We orient each $C_i$ to get a directed cycle. Now for each singleton edge $e\in E(G)$
such that $e$ is contained in some $C_i$ for $i\in [1,m]$, we orient $e$ and its duplication  $e'$ in $G^*$
the same direction as it is on $C_i$. 
For each singleton edge  $e\in E(G)$ such that $e$ is contained in some $P_i$
for $i\in [1,\ell]$, we direct $e$ and its duplication $e'$ in $G^*$
in opposite directions. Now we have oriented all the edges of $G^*$. Call the resulting directed graph $D$. 
By the construction, such an  orientation satisfies the two properties above.

We then construct   a bipartite graph $H$ associated with $D$. 
Let   $V(H)=\{v^+, v^-: v\in V(D)\}$.  An arc $uv$ in $D$ is corresponding to an edge joining $u^+$ and $v^-$ in $H$. 
By this construction, for any vertex $v\in V(D)$, we have 
\begin{enumerate}[(1)]
	\item  $d_H(v^+)=d_D^+(v)=d_H(v^-)=d_D^-(v)=rs$ for any $v\in V(D)$;  
	\item $\mu(H) \le r$; 
\end{enumerate}

Let  $f(v)=s$ and $g(uv)=1$
for all  $u, v\in V(H)$.   As $H$ is $r s$-regular with $\mu(H) =r$, 
we find an $fg$-coloring of $H$ using $r$ colors by Lemma~\ref{lem:fg-b}. 
In each of the color classes,  identifying $v^+$ and $v^-$ for each $v\in V(G^*)$   
gives a decomposition of $G^*$ into $r$ spanning subgraphs $H_1, \ldots, H_r$, 
where each $H_i$ is $2s$-regular with maximum multiplicity 2 (since the maximum multiplicity of $G^*$ is $2r$ and each color classes of $H$ is a simple graph).  

Let $T^*=P_1\cup \ldots \cup P_\ell$.  As $T^*$ 
contains at most $n-1$ edges, at least one of those subgraphs $H_1, \ldots, H_r$ has fewer than  $n/3$ edges of  $T^*$. 
Relabeling the $H_i$, if necessary, we can therefore assume  $H_1$  contains fewer than  $n/3$ edges of  $T^*$.  
Thus for every $v\in V(H_1)$, 
$v$ is incident in $H_1$ with at most $n/3$ edges from $E(H_1)\cap E(T^*)$. 



Ideally we would now like to partition the graph  $H_{1}$  into two spanning simple graphs, each with maximum degree $s$.  However, that may not be possible; for example if $H_{1}$ has every edge having multiplicity two, except for 3 singleton edges that form a triangle using vertices of degree $2s$, this is impossible.  However we claim that we can find a decomposition of $ H_{1}$  into a Hamilton cycle (and thus two perfect matchings) and two spanning simple graphs, each with maximum degree $s-1$.  For every singleton edge $e\in E(G)$, if both $e$ and its duplication $e'$ 
are contained in $H_{1} $, we delete both $e$ and $e'$ from $H_{1}$.  
Denote the resulting graph by $H'_1$. 

\begin{CLA}\label{claim:multiple-edges}
For each $v\in V(H_1')$,  $v$ is incident in $H_1'$ to either  at least $n/3+2$ singleton edges or at least $16$
multiple edges. 
\end{CLA}

\pf Let $v\in V(H_1')$ be any vertex.   We suppose $v$ is adjacent in $H_1'$ to at most  $15$
multiple edges.  Then as $\delta(H_1) \ge n+32$, 
we then know that 
 $v$ is incident in $H'_1$ to at least  $n+32-2(n/3+15)=n/3+2$ singleton edges 
 of $H_1'$. 
\qed 

%

\begin{CLA}\label{claim:H1-prime-con}
	The graph $H_1'$ is connected. 
\end{CLA}

\pf In fact, we show that $H_1'$ has a Hamilton cycle. 
Let $J$  and $J^*$ be the underlying simple graphs of $H'_1$  and $H_1$, respectively. 
Assume  $d_J(v_1) \le \ldots \le d_J(v_{n})$, where $\{v_1,\ldots, v_{n}\}=V(H'_1)$.  Suppose to the contrary that $J$ does not have a Hamilton cycle. 
Then by Theorem~\ref{lem:chvatal's-theorem}, there exists $i<[1,n/2-1]$ such that $d_J(v_i) \le i$ and $d_J(v_{n-i}) \le n-i-1$. 
By Claim~\ref{claim:multiple-edges}, we know that in $H'_1$, each vertex is incident with at least  $16$ distinct vertices. 
Thus $\delta(J) \ge 16$ and so $i \ge 17$.   Recall that  $T^*$ 
is a forest of $G$.  Thus 
\begin{eqnarray*}
	d_{J}(v_i) &\ge&  \frac{1}{i} \sum\limits_{j=1}^id_J(v_j) \ge 
	\frac{1}{i}\left (\sum\limits_{j=1}^i d_{J^*}(v_i) -(2n-2) \right)  \\
	&\ge  & (n/2+15)-\frac{1}{i}(2n-2). 
\end{eqnarray*}
This implies $d_{J}(v_i) \ge \frac{3}{8}n$  as $i\ge 17$. Thus $i>\frac{3}{8}n$. Now by the same argument above we get 
$$
d_{J}(v_i)>(n/2+16)-\frac{1}{i}(2n-2)>n/2, 
$$
a contradiction. 
\qed 

Since $H_1$ is even and we only deleted some doubleton edges of $H_1$ to get $H_1'$, 
we know that $H_1'$ is even. Thus the subgraph, call it  $R_1'$,  of $H_1'$ induced on its singleton edges 
is even. 
By applying Lemma~\ref{lem:spanning forest} to each component of $R_1'$, we can find a 
forest  $T$    of  $H'_1$ that is formed by its singleton edges
and has as many edges as possible so  $d_T(v) \le 1+\frac{1}{2}d_{R_1'}(v) \le n/6+2$ for each $v\in V(R_1')$ if $d_{R_1'}(v) \ge n/3+2$.
Thus $d_{R_1'-E(T)}(v) \ge n/6\ge 16 $ by $n\ge 96$. 
Let $H_1^*=H'_1-E(T)$.   

\begin{CLA}\label{claim:h-cycle}
The graph $ H_1^*$  has a Hamilton cycle.
\end{CLA}

\pf Let $J$  and $J^*$ be the underlying simple graphs of $H'_1$  and $H_1$, respectively. 
Assume  $d_J(v_1) \le \ldots \le d_J(v_{n})$, where $\{v_1,\ldots, v_{n}\}=V(H^*_1)$.  Suppose to the contrary that $J$ does not have a Hamilton cycle. 
Then by Theorem~\ref{lem:chvatal's-theorem}, there exists $i<[1,n/2-1]$ such that $d_J(v_i) \le i$ and $d_J(v_{n-i}) \le n-i-1$. 
Note that  $\delta(J) \ge 16$ by Claim~\ref{claim:multiple-edges} and the fact that $d_{R_1'-E(T)}(v) \ge n/6\ge 16 $ if $d_{R_1'}(v) \ge n/3+2$. 
Thus  $i \ge 17$.   Recall that the union  $T^*$ of $P_1, \ldots, P_\ell$ 
is a forest of $G$.   Then
\begin{eqnarray*}
	d_{J}(v_i) &\ge&  \frac{1}{i} \sum\limits_{j=1}^id_J(v_j) \ge 
	\frac{1}{i}\left (\sum\limits_{j=1}^i d_{J^*}(v_i) -(4n-4) \right)  \\
	&\ge  & (n/2+16)-\frac{1}{i}(4n-4). 
\end{eqnarray*}
This implies $i>n/4$. Now by the same argument above we get 
$$
d_{J}(v_i)>(n/2+16)-\frac{1}{i}(4n-2)>n/2, 
$$
a contradiction. 
\qed 

Thus $H_1^*$  has a Hamilton cycle, call it $C$.  Now consider $H_{1} - E(C)$.  Partition its edges to form two simple graphs $H_{11}$ and $H_{12}$ as follows.  First place one of each of  the doubleton edges of $H_{1} - E(C) $ into these two graphs.  
Let $R$ be the subgraph of $H_1-C$ formed by all its singleton edges. Note that by our choice of $C$, the graph $R$ is connected. 
This is because (1)
$E(T) \subseteq E(R)$;  and (2) the edges of $H_1$ joining components of $T$ (if  $T$ has more than one components) are all doubleton edges, 
thus in this case the cycle $C$ contains edges of $H_1$ that are joining components of $T$ and so  deleting edges in $C$ leaves singleton edges in  $ H_{1} - E(T)$ that are still joining the components of $T$. 
Then note that all the vertices of $R$ are of even degree. This is because $H_1$ is $2s$-regular, $C$
is 2-regular, and we only removed doubleton edges from $H_1-E(C)$ to get $R$. 
As $R$ is connected with all vertices of even degree, $R$ has an Euler tour by Theorem~\ref{Euler}. 
Since $n$ is even and $2s$ is even, we know that $e(H_1)$ is even and so 
$R$ has an even number of edges.  
Now placing the edges of  $R$ from one of its Euler tour alternately in $H_{11}$ and $H_{12}$ gives the desired decomposition, each graph is $(s-1)$-regular.

We now let $G_1=G-E(C)-E(H_{11})$. As  $\mu(G^*)=2r$, we know that  each $H_i$ contains two edges between any two vertices $u,v\in V(G^*)$
with $e_{G^*}(u,v)=2r$. 
As $H_{12}$ is a simple graph, we know that $C\cup H_{11}$ contains at least one edge  between any two vertices $u,v\in V(G^*)$
with $e_{G^*}(u,v)=2r$.  Thus $C\cup H_{11}$ contains at least one edge  between any two vertices $u,v\in V(G)$
with $e_{G}(u,v)=r$.  
Thus $\mu(G_1) \le r-1$. As $C\cup H_{11}$ is spanning and  $(s+1)$-regular, we know 
that $G_1$ is $((r-1)s-1)$-regular. By Dirac's Theorem, the underlying simple graph of $G_1$ contains $r-2$ edge-disjoint perfect matchings. 
Delete those  $r-2$  perfect matchings from $G_1$ to get $G_2$. Then $G_2$
is $(r-1)(s-1)$-regular. Applying Lemma~\ref{lem:decomposition}, $G_2$
can be decomposed into $(r-1)/2$ edge-disjoint Hamilton cycles and $(r-1)$  spanning $(s-2)$-regular simple graphs.  
Thus 
	$G$ can be decomposed into $r-2$ perfect matchings,    $(r+1)/2$ 
Hamilton cycles, one  spanning $(s-1)$-regular simple graph, and $(r-1)$  spanning  $(s-2)$-regular simple graphs.  
\qed 

\proof[Proof of Theorem~\ref{thm:1-factorization}]
We choose $n_0$ to be the same as that stated in Theorem~\ref{thm:1-factorization-proof}. Let $G$ be a  $k$-regular graph 
with maximum multiplicity $r$ and $k\ge r(n/2+18)$. Applying Theorem~\ref{thm:Dirac}
to the underlying simple graph of $G$, if necessary, by removing  at most $r-1$ edge-disjoint perfect matchings from $G$, 
we may assume  that $k=rs$ for some integer $s\ge n/2+17$. 

If $r$ is even, applying Lemma~\ref{lem:decomposition}, we decompose $G$
into $r/2$ edge-disjoint Hamilton cycles and $r$ spanning $(s-1)$-regular simple graphs.
Each of the Hamilton cycle can be decomposed into $2$ perfect matchings and each of the $r$
simple graphs has a 1-factorization by Theorem~\ref{thm:1-factorization-proof}. 
Thus $G$ has a 1-factorization into $r+r(s-1)=rs$ 1-factors. 

Now suppose that $r$ is odd. If $s$ is even, by Lemma~\ref{lem:decompose2}(a), 
$G$ can be decomposed into $r-2$ perfect matchings,    $(r+1)/2$ 
Hamilton cycles, one spanning  $(s-1)$-regular simple graph, and $(r-1)$  spanning $(s-2)$-regular simple graphs.  
 If $s$ is odd,  by Lemma~\ref{lem:decompose2}(b),  
$G$ can be decomposed into $2r-2$ perfect matchings,    $(r+1)/2$ 
Hamilton cycles, one  spanning $(s-2)$-regular simple graph, and $(r-1)$  spanning  $(s-3)$-regular simple graphs.  
Again as each of the Hamilton cycle can be decomposed into $2$ perfect matchings and each of the  
simple graphs has a 1-factorization by Theorem~\ref{thm:1-factorization-proof}, we know that 
 $G$ has a 1-factorization. 
\qed 

%
%
%

 \proof[Proof of Theorem~\ref{thm:overfull-present}]  
 We choose $n_0$ to be the same as that stated in Theorem~\ref{thm:1-factorization-proof}. 
 We assume  that  $G$  is not regular, else the result follows from Theorem~\ref{thm:1-factorization}.  By Lemma~\ref{lemma:overfull-subgraph}, $G$  contains a vertex  $v$  such that  $\omega(G- v) = \omega(G)$. Let 
 $$
 \omega(G -v) =\delta(G)+k+t,
 $$
 where $k\ge 0$ is an integer  and  $0 \le t<1$ is rational of the form $\frac{b}{n-2}$, where $0\le b<n-2$ is even. 
 Let  $w$  be a vertex with minimum degree among all vertices other than  $v$, which can be of degree $\Delta(G)$ 
 if $v$ is the only vertex of degree less than $\Delta(G)$ in $G$.  
   By Dirac's Theorem,  $G - \{v, w\}$  has a perfect matching; let  $M$  be a subset of such a matching with exactly  $t\cdot (n-2)/2 $ edges, and consider $G_0=G-M$. 
 The  graph  $G_0$ has  $G_0- v$  as a $\Delta(G_0)$-full/overfull subgraph, and the number of edges of  $G_0-v$  is now a multiple of $(n-2)/2$, specifically $(\delta(G_0)+k)(n-2)/2$.   Let  $w_0$  be a vertex of minimum degree in  $G_0$  other than $v$.  As before,  $G_0-\{v,w_0\}$  has a perfect matching  $M_0$; we remove that matching from  $G_0$.  As we continue this process, adjusting the choice of  $w_0$  each time as necessary, by Lemmas~\ref{lemma:overfull-subgraph2} and~\ref{lemma:overfull-subgraph},  we have at most two vertices of minimum degree at each step, one of which must be vertex  $v$  whose degree does not change throughout the procedure. 
 After the removal of  $k$  such perfect matchings the graph  $G^*$  obtained has minimum degree $\delta(G)$, and  $G^*-v$  has exactly    $\delta(G)(n-2)/2$ edges.  It is straightforward then to see that  $G^*$  must be   $\delta(G)$-regular.  This graph has a 1-factorization by Theorem~\ref{thm:1-factorization}, and combining this 1-factorization with the previously removed matchings gives a decomposition of  $G$  into  matchings.  The result follows as $\chi'(G) \ge \omega(G)$. 
 
 The consequence part  of the statement is clear as when $G$ contains a $\Delta(G)$-full subgraph, we have $\lceil \omega(G) \rceil =\Delta(G)$. 
  \qed

 \section{Proof of Theorem~\ref{thm:1}}

\subsection{Special cases}

\begin{LEM}\label{lema:small-deficiency}
 There exists an $n_0\in \mathbb{N}$ such that the following holds. Let $n, r \in \mathbb{N}$
be such that $n\ge n_0$ is even and $r=\mu(G)$.   
 If $G$ has a minimum degree vertex $v^{*}$ such that $e(G- v^{*}) =  \Delta(G)(n-2) /2-s$
 for some $s\in [0,6r]$ and $\delta(G) \ge r(n/2+7)+s$, then $\chi'(G) = \Delta(G) $. 
 \end{LEM}

\pf   If $ e(G- v^{*}) =   \Delta(G) (n-2)/2$, then  $G - v^{*}$ is $ \Delta(G)$-full in $G$ and the result follows by  Theorem~\ref{thm:overfull-present}.  So we assume that $ e(G - v^{*})= \Delta(G)(n-2) /2- s$  for  $s\in[1,6r]$. 
Since  
\begin{eqnarray*}
\sum_{v\in V(G)}(\Delta(G)-d_G(v))&=&\sum_{v\in V(G-v^*)}(\Delta(G)-d_{G^*}(v))-d_G(v^*)+(\Delta(G)-d_{G}(v^*))\\
&=&(n-1)\Delta(G)-2e(G-v^*)-d_G(v^*)+(\Delta(G)-d_{G}(v^*))\\ 
&=&\Delta(G)+2s-d_G(v^*)+(\Delta(G)-d_{G}(v^*)) \\
	 &=& 2(\Delta(G)-d_{G}(v^*))+2s, 
\end{eqnarray*}
and $\Delta(G)-d_{G}(v^*)=\max\{\Delta(G)-d_{G}(v): v\in V(G)\}$, we know that 
 there exist  at least two other vertices $u,w$ in $G$ with degree less than $\Delta(G)$.  By Dirac's Theorem applied to the simple graph underlying $G$, there is a Hamilton cycle, and therefore a perfect matching $M$, in $G - u - w$.  We remove $M$ from $G$ and let $G_1=G-M$. 
 
Note that $\delta(G_1)=\delta(G)-1=d_{G_1}(v^*)$. 
   Furthermore, $e(G_1-v^*)=\Delta(G)(n-2)/2 - s-(n-4)/2=(\Delta(G)-1)(n-2)/2-(s-1)=\Delta(G_1)(n-2)/2-(s-1) \le \Delta(G_1)(n-2)/2$.  
    Since $\delta(G_1)=\delta(G)-1\ge r(n/2+7)+s-1$,
 we can repeat the process with $G_1$ in the place of $G$. 
   In general,  we repeat this process $s$ times  to yield a graph  $G_s$ of maximum degree $\Delta(G) - s$  in which  $ G_s - v^{*} $ is $\Delta(G_s)$-full with $\delta(G_s) \ge r(n/2+7)$. By Theorem~\ref{thm:overfull-present},
   we have $\chi'(G_s)=\Delta(G)-s$. As each time a matching was only removed 
   when we repeated the process, we then know that   $\chi'(G) = \Delta(G)$.  
   \qed


\subsection{Proof of Theorem~\ref{thm:1}}

\pf  We choose $n_0$ to be at least the maximum of the $n_0$ stated  in Theorem~\ref{thm:1-factorization-proof} and
the $n_0$ with respect to $0.9\ve$ as stated in Theorem~\ref{thm:plantholt-shan}, and such that $1/n_0\ll \ve$.  
Let  $G$ be a graph of  even order  $n \ge  n_{0}$,   maximum multiplicity  $r$, minimum degree  $\delta > r(1+\ve)n/2$, maximum degree $\Delta $.  Denote by $V_\delta$ the set of minimum degree vertices of $G$. 

By Lemma~\ref{lema:small-deficiency},  we assume 
\begin{equation}\label{eqn:e-of-g-v1}
	e(G - v) < \Delta(n-2)/2 - 6r \quad \text{for any $v\in V(G)$.}
\end{equation} 
Therefore, if two vertices with degree less than $\Delta$  are not adjacent in  $G$ with $r$ edges between them, we may add an edge between them without creating an overfull subgraph, or increasing $\Delta$. We iterate this edge-addition procedure.  If at some point we create a  $\Delta$-full subgraph, the result follows by Theorem~\ref{thm:overfull-present}.  Otherwise, we reach a point where we
may now assume that in  $G$  all vertices with degree less than  $\Delta$ are mutually adjacent with $r$ edges between them. 
Thus by Lemma~\ref{lem:numberofD}, we know that $|V_\Delta| \ge n/2+1$. 

We may assume that the maximum degree $\Delta$ is a multiple of $r$,  say  
$$\Delta = rk \quad \text{for some positive integer $k$, and $k$ is even when $r$ is odd}.$$  
For otherwise,  when $r$ is even, we can remove at most $r-1$ edge-disjoint perfect matchings to achieve 
the property;  and 
when $r$ is odd, 
we can remove at most $2r-1$  edge-disjoint perfect matchings to achieve that property. 
The existence of perfect matchings is guaranteed as the underlying simple graph of $G$ is hamiltonian by Dirac's Theorem. 
Note that deleting a perfect matching maintains the 
inequality in~\eqref{eqn:e-of-g-v1} with respect to the 
resulting graph.   The minimum degree  $\delta$
 satisfies $\delta \ge r(1+ \ve)n/2-2r \ge  r(1+ \ve')n/2$, where $\ve'=0.9\ve$.

We first form a supergraph $ G^{*}$  that contains $G$ as a subgraph as follows: 
\begin{enumerate}[---]
	\item Add two vertices $z_1$ and $w_1$ to $G$.
	\item For any vertex  $u$ of $G$ whose degree is not a multiple of $r$, add just enough parallel $uz_1$ edges so the degree of $u$ is a multiple of $r$ in $G^{*}$. 
	\item Finally, if necessary, add enough parallel edges (at most $r-1$ edges) between $ z_1$ and $w_1$ so that  the degree of $z_1$ in $G^{*}$ is also a multiple of $r$.
\end{enumerate}

Since $|V_\Delta| \ge n/2+1$, we have 
$$d_{G^*}(z_1) \le (r-1)(n/2-1)+r-1 \le (r-1)n/2 <\delta.$$ 
We now separate the proof into two cases according to whether  or not $r$ is even. 

\smallskip 

{\bf \noindent Case 1: $r$ is even}. 

\smallskip

We find an $fg$-coloring of $G^{*}$ by letting 
\begin{numcases}{f(v)=}
\frac{d_{G^*}(v)}{r/2} & \text{if $v\in V(G^*)\setminus\{w_1\}$;}  \nonumber \\
 2 &  \text{if $v=w_1$;} \nonumber 
\end{numcases}
and 
$$
g(uv)=2 \quad \text{for any  pair $u,v\in V(G^*)$}. 
$$
 By Lemma~\ref{lem:fg}, since each $f, g $ value is even, we can obtain such a coloring using $\Delta/(\Delta / \frac {r} {2}) =  r/2 $ colors. This gives us a partition of $G^{*}$  into $ r/2 $ subgraphs $H_{1}^{*}, H_{2}^{*}, ..., H_{r /2}^{*}$, each with maximum multiplicity at most 2, and maximum degree $2k$. Let  $H_i=H_i^*-\{z_1,w_1\}$ for each $i\in [1,r/2]$.

 We first claim that none of these graphs $H_i$ has a $2k$-overfull subgraph.  To see this,  note that $e(G) = e(G - v^{*}) + \delta  \leq \Delta(n-2)/2 - 6r + \delta$ for any $v^*\in V_\delta$  by~\eqref{eqn:e-of-g-v1}. 
  Because in our decomposition each vertex  $v$ of $H_{i}^{*} $ other than $w_1$ has degree $f(v)$, and $w_1$ has degree at most 2, it follows that for any $i,j$, the edge cardinalities of  $H_{i}$ and $ H_{j}  $  differ by at most 1.  Thus each $H_{i}$ has at most $ \left \lceil   ( e(G - v^{*} ) + \delta)/ (r /2) \right \rceil$ edges. In addition, since in the coloring $w_1$ has degree at most 2 in each color class, the degree of any vertex $v$ can differ by at most 2 in any two subgraphs $H_{i}, H_{j}$.  It follows then that for any $H_{i}$, the number of edges in any vertex-deleted subgraph is at most  
  \begin{eqnarray*}
   &&\left \lceil   ( e(G - v^{*}) + \delta)/ (r /2) \right \rceil - (\delta/(r /2)) +2  \\
    &\leq&  (kr (n-2)/2- 6r )/ (r /2) + 2\delta/r +1 - 2\delta/r + 2 \leq 2k(n-2)/2 - 9. 
  \end{eqnarray*}
  Thus $H_i$ contains no $2k$-overfull subgraph by Lemma~\ref{lemma:overfull-subgraph}.   Each $H_i$ satisfies $\delta(H_i) \ge (1+\ve' ) n$.

Just as in the proof of Lemma~\ref{lem:decompose2}, we find a decomposition of $ H_{i}$  into a Hamilton cycle (and thus two perfect matchings) and two simple graphs, each with maximum degree $k-1$.  
To get this partition, first let  $T_{i}$ be a spanning forest  of the subgraph of $H_{i}$ that is formed by its singleton edges
and has as many edges as possible.  We show below that $ H_{i} - E(T_{i})$  has a Hamilton cycle. 

\begin{CLA}\label{claim:h-cycle2}
Each graph $ H_{i} - E(T_{i})$  has a Hamilton cycle.
\end{CLA}

\pf Let  $J^*$ and $J$ be the underlying simple graph of  $H_i$ and $H_i-E(T_i)$, respectively.  
Assume   $d_J(v_1) \le \ldots \le d_J(v_{n})$, where $\{v_1,\ldots, v_{n}\}=V(H_i)$.  Suppose to the contrary that $J$ does not have a Hamilton cycle. 
Then by Theorem~\ref{lem:chvatal's-theorem}, there exists $i<[1,n/2-1]$ such that $d_J(v_i) \le i$ and $d_J(v_{n-i}) \le n-i-1$. 
Since $\delta(H_i) \ge (1+\ve')n$, we know that in $H_i$, each vertex is incident with at least $\ve' n$ multiple edges. 
Thus $\delta(J) \ge \ve' n$ and so $i >\ve' n$. Then as $\sum\limits_{j=1}^{n} d_{T_i}(v_j) \le 2n-2$, we get 
\begin{eqnarray*}
d_{J}(v_i) &\ge&  \frac{1}{i} \sum\limits_{j=1}^id_J(v_j) \ge 
\frac{1}{i}\left (\sum\limits_{j=1}^i d_{J^*}(v_i) -(2n-2) \right)  \\
 &\ge  & ( 1+ \ve' )n/2-\frac{1}{i}(2n-2) >n/2, 
 \end{eqnarray*}
showing a contradiction to  $d_J(v_i) \le i <n/2$.
\qed 

Thus $ H_{i} - E(T_{i})$  has a Hamilton cycle, call it $C_i$.  Now consider $H_{i} - E(C_i)$.  Partition its edges to form two simple graphs $H_{i1}$ and $H_{i2}$ as follows.  First place one of each of the the doubleton edges of $H_{i} - E(C_i) $ into these two graphs.  
Let $R_i$ be the subgraph of $H_i-C_i$ formed by all its singleton edges. Note that by our choice of $C_i$, the graph $R_i$ is connected. 
This is because (1)
$E(T_i) \subseteq E(R_i)$;  and (2) the edges of $H_i$ joining components of $T_i$ (if $T_i$ has  more than one components) are all doubleton edges, 
thus in this case the cycle $C_i$ contains edges of $H_i$ that are joining components of $T_i$ and so  deleting edges in $C_i$ leaves singleton edges in  $ H_{i} - E(T_{i})$ that are still joining the components of $T_{i}$. 
We add a new vertex $y$ to $R_i$ and add an edge between $y$ and each odd degree vertex of $R_i$, 
 calling the new graph $R^{*}_i$.  Then  we find an Euler tour of  $R^{*}_i$ that begins and ends at a vertex that is not of maximum degree in the graph $H_i$ ($H_i$ is not regular as any vertex-deleted subgraph of it has at most $2k(n-2)/2-9$  edges).
 Placing the edges of  $R_i$ from this tour alternately in $H_{i1}$ and $H_{i2}$ gives the desired decomposition, each graph having maximum degree $k-1$.

Finally we show that neither $H_{i1}$ nor $H_{i2}$ has a $(k-1)$-overfull subgraph.  Recall from before that $H_{i}$ has at most $2k(n-2)/2 - 9$ edges in any vertex-deleted subgraph.  Therefore  $H_{i} - E(C_i)$ has at most $ (2k-2)(n-2) /2- 9$ edges in any vertex-deleted subgraph. By construction, the graphs  $H_{i1}$ and $H_{i2}$ differ in size by at most 1, and in degree at any vertex by at most 2. Thus, arguing as before, the number of edges in any vertex-deleted subgraph of  $H_{i1}$ or $H_{i2}$  is at most  $(k-1)(n-2)/2-9 + 3 <(k-1)(n-2) /2$, so that  $H_{i1}$ and $H_{i2}$  contain no $(k-1)$-overfull subgraph.  It follows by Theorem~\ref{thm:plantholt-shan} that  $ \chi'(H_{i1}) =   \chi'(H_{i2}) = k-1$.  Consequently  $\chi'(G) = \Delta(G)$, and the result follows.

\smallskip 

{\bf \noindent Case 2: $r$ is odd}. 

\smallskip

We first form a supergraph $ G^{**}$  based on $G^*$:  
\begin{enumerate}[---]
	\item Add two vertices $z_2$ and $w_2$ to $G^*$.
	\item For any vertex  $u$ of $G^*$ with $u\ne w_1$ whose degree is not an even multiple of $r$ (so  the degree of $u$ in $G^*$ is at most $\Delta-r$), add exactly $r$ parallel $uz_2$ edges so the degree of $u$ is an even multiple of $r$ in $G^{**}$. 
	\item Finally, if necessary, add $r$ parallel edges between $ z_2$ and $w_2$ so that  the degree of $z_2$ in $G^{**}$ is also an even multiple of $r$.
\end{enumerate}

Since $|V_\Delta| \ge n/2+1$, we have 
$$d_{G^{**}}(z_2) \le  r(n/2-1)+r = r n/2 <\delta,$$ 
 $d_{G^{**}}(w_1) \le r-1$ and $d_{G^{**}}(w_2) \le r$.

We double all the edges of $G^{**}$ to get $G'$. Then $G'$ has  maximum multiplicity $2r$. By the same argument 
as in the proof of Lemma~\ref{lem:decompose2}, 
we know  that $G'$ has an orientation $D$
such that 
\begin{enumerate}[(1)]
	\item  $d_D^+(v)=d_D^-(v)$ for any $v\in V(D)$;  
	\item for any $u,v\in V(D)$, there are at most $r$ arcs of $D$ from $u$ to $v$; 
	\item Let $E=\{e\in A(D): \text{ $e$ is singleton in $G$ and $e$ and its duplication $e'$ form a 2-cycle  in $D$}\}$, where $A(D)$ is the set of arcs of $D$.  Then $D[E]$ 
	is a forest. 
\end{enumerate}

We then construct   a bipartite graph $H$ associated with $D$. 
Let   $V(H)=\{v^+, v^-: v\in V(D)\}$.  An arc $uv$ in $D$ is corresponding to an edge joining $u^+$ and $v^-$ in $H$. 
By this construction, for any vertex $v\in V(D)$, we have 
\begin{enumerate}[(1)]
	\item  $d_H(v^+)=d_D^+(v)=d_H(v^-)=d_D^-(v)$ for any $v\in V(D)$;  
	\item $\mu(H) \le r$; 
	\item The union of the  edges corresponding to  $E$ form a forest  $T^*$ in $H$. 
\end{enumerate}

We find an $fg$-coloring of $H$ by letting 
\begin{numcases}{f(v)=}
	\frac{d_{H}(v)}{r} & \text{if $v\in V(H)\setminus\{w_1^+, w_1^-, w_2^+, w_2^-\}$;}  \nonumber \\
	1 &  \text{if $v\in \{w_1^+, w_1^-, w_2^+, w_2^-\}$;} \nonumber 
\end{numcases}
and 
$$
g(uv)=1 \quad \text{for any  pair $u,v\in V(H)$}. 
$$

 As $d_{H}(v)$ is a multiple of  $r$ for any $v\in V(H)\setminus\{w_1^+, w_1^-, w_2^+, w_2^-\}$,  $d_H(v) \le r$
 for any $v\in \{w_1^+, w_1^-, w_2^+, w_2^-\}$, and $\mu(H)=r$,
 we find an $fg$-coloring of $H$ using $r$ colors by Lemma~\ref{lem:fg-b}. 
In each of the color classes,  identifying $v^+$ and $v^-$ for each $v\in V(G^{**})$   
gives a decomposition of $G^{**}$ into $r$ spanning subgraphs $H^*_1, \ldots, H^*_r$, 
where each $H^*_i$ has maximum degree $2k$
and maximum multiplicity 2 (since the maximum multiplicity of $G^{**}$ is $2r$ and each color class of $H$ is a simple graph).  
Let $H_i=H_i^*-\{z_1,w_1,z_2,w_2\}$ for each $i\in [1,r]$.

  We first claim that none of these graphs  $H_i$ has a $2k$-overfull subgraph.  To see this,  note that $e(G) = e(G - v^{*}) + \delta  \leq \Delta(n-2)/2 - 6r + \delta$ for any $v^*\in V_\delta$ by~\eqref{eqn:e-of-g-v1}.  
 Because in our decomposition each vertex  $v$ of $H_{i}^{*} $ other than $w_1$ and $w_2$ has degree $f(v)$, and $w_1, w_2$ have degree at most 2, it follows that for any $i,j$, the edge cardinalities of  $H_{i}$ and $ H_{j} $  differ by at most 2.  Thus each $H_{i}$ has at most $   ( e(G - v^{*} ) + \delta)/ r   +2$ edges. In addition, since in the coloring $w_1,w_2$ have degree at most 2 in each color class, the degree of any vertex $v$ can differ by at most 4 in any two subgraphs $H_{i}, H_{j}$.  It follows then that for any $H_{i}$, the number of edges in any vertex-deleted subgraph is at most  
 \begin{eqnarray}
 	&&   2( e(G - v^{*}) + \delta)/ r  +2- (\delta/r) +4  \nonumber \\
 	&\leq&  2(rk (n-2)/2 - 6r )/ r  + \delta/r +2 - \delta/r + 4 \leq 2k(n-2)/2 - 6.  \label{eqn:Hi-v}
 \end{eqnarray}
 Thus $H_i$ contains no $2k$-overfull subgraph by Lemma~\ref{lemma:overfull-subgraph}.  Each $H_i$ satisfies $\delta(H_i) \ge (1+\ve' ) n$. 
 
 By the same argument as in the proof of Lemma~\ref{lem:decompose2}, 
 for every $v\in V(H_1)$,  we may assume that 
 $v$ is incident in $H_1$ with at most $n/3$ edges from $E(H_1)\cap E(T^*)$. 
   We claim that we can find a decomposition of $ H_{1}$  into a Hamilton cycle (and thus two perfect matchings) and two spanning simple graphs, each with maximum degree $s-1$.  For every singleton edge $e\in E(G)$, if both $e$ and its duplication $e'$ 
 are contained in $H_{1} $, we delete both $e$ and $e'$ from $H_{1}$.  
 Denote the resulting graph by $H'_1$.  
 
 \begin{CLA}\label{claim:multiple-edges2}
 	For each $v\in V(H_1')$,  $v$ is incident in $H_1'$ to either  at least $n/3+2$ singleton edges or at least $ \lfloor \ve' n \rfloor$
 	multiple edges. 
 \end{CLA}
 
 \pf Let $v\in V(H_1')$ be any vertex.   We suppose $v$ is adjacent in $H_1'$ to at most  $\ve' n-1$
 multiple edges.  Then as $\delta(H_1) \ge n+2\ve' n$, 
 we then know that 
 $v$ is incident in $H'_1$ to at least  $n+2\ve' n-2(n/3+ \ve' n-1)=n/3+2$ singleton edges 
 of $H_1'$. 
 \qed 
 
 %
 
 \begin{CLA}\label{claim:H1-prime-con2}
 	The graph $H_1'$ is connected.   
 \end{CLA}

 \pf  In fact, we show that $H_1'$ has a Hamilton cycle. 
 Let $J$  and $J^*$ be the underlying simple graphs of $H'_1$  and $H_1$, respectively. 
 Assume  $d_J(v_1) \le \ldots \le d_J(v_{n})$, where $\{v_1,\ldots, v_{n}\}=V(H'_1)$.  Suppose to the contrary that $J$ does not have a Hamilton cycle. 
 Then by Theorem~\ref{lem:chvatal's-theorem}, there exists $i<[1,n/2-1]$ such that $d_J(v_i) \le i$ and $d_J(v_{n-i}) \le n-i-1$. 
 By Claim~\ref{claim:multiple-edges2}, we know that in $H'_1$, each vertex is incident with at least  $ \lfloor \ve' n \rfloor$ distinct vertices. 
 Thus $\delta(J) \ge  \lfloor \ve' n \rfloor$ and so $i \ge  \lfloor \ve' n \rfloor+1$.   Recall that  $T^*$ 
 is a forest of $G$.  Thus 
 \begin{eqnarray*}
 	d_{J}(v_i) &\ge&  \frac{1}{i} \sum\limits_{j=1}^id_J(v_j) \ge 
 	\frac{1}{i}\left (\sum\limits_{j=1}^i d_{J^*}(v_i) -(2n-2) \right)  \\
 	&\ge  & (n/2+\ve' n)-\frac{1}{i}(2n-2)>n/2, 
 \end{eqnarray*}
  a contradiction. 
 \qed

 
 Let $R'$ be the subgraph of $H_1'$ induced on its singleton edges. 
 Since $H_1$ contains doubleton edges and singleton edges only, 
 the parities of the degrees of each vertex in $H_1$, $H_1'$,  and in $R'$ are the same.  

If $R'$ is even,  applying Lemma~\ref{lem:spanning forest} to each component of $R'$, we can find a 
 forest  $T$    of  $H'_1$ that is formed by its singleton edges
 and has as many edges as possible  such that   $d_{R'-E(T)}(v)  \ge \ve' n$ 
 if $d_{R'}(v) \ge n/3+2$. 
 
If $R'$ is not even, then we add a new vertex $w$ to $R'$ and join an edge between $w$
and each vertex of $R'$ that is of odd degree. Denote the resulting graph by $R^{''}$. 
Now applying Lemma~\ref{lem:spanning forest} to each component of $R^{''}$, we can find a 
forest  $T$    of  $R^{''}$ with as many edges as possible  such that  $d_{R^{''}-E(T)}(v) \ge \ve' n$
if $d_{R^{''}}(v) \ge n/3+2$ for $v\in V(H_1')$.  By adding edges of $R^{''}$ incident to $w$
and deleting other edges if necessary (this will maintain the property that $d_{R^{''}-E(T)}(v) \ge \ve' n$
if $d_{R^{''}}(v) \ge n/3+2$ for $v\in V(H_1')$), we may assume that all edges incident with $w$ in $R^{''}$
are contained in $T$. 

Let $H_1^*=H'_1-E(T)$.   
By a similar argument as in Case 1 but using $4n-4$ in the place of $2n-2$,  we know that $H_1^*$  has a Hamilton cycle, call it $C$.  Now consider $H_{1} - E(C)$.  Partition its edges to form two simple graphs $H_{11}$ and $H_{12}$ as follows.  First place one of each of  the doubleton edges of $H_{1} - E(C) $ into these two graphs.  

 Let $R$ be the subgraph of $H_1-C$ formed by all its singleton edges together 
 with the edges of $R^{''}$ that are incident with $w$ if they exist. 
  Note that by our choice of $C$, the graph $R$ is connected. 
 This is because (1)
 $E(T) \subseteq E(R)$;  and (2) the edges of $H_1$ joining components of $T$ (if  $T$ has more than one components) are all doubleton edges (recall that $T$ contains all edges incident with $w$ in $R^{''}$), 
 thus in this case the cycle $C$ contains edges of $H_1$ that are joining components of $T$ and so  deleting edges in $C$ leaves singleton edges in  $ H_{1} - E(T)$ that are still joining the components of $T$. 
 
 We then claim that $R$ is even.   Let $R^*=R'$ if $R'$ is even 
 and $R^*=R^{''}$ otherwise. 
 Note that at each vertex $v\in V(R)\cap V(H_1^*)$, there are three possibilities of constitution of the  two edges of $C$ that
 are  incident with $v$: both are two doubleton edges of $H_1^*$ ($d_R(v)=d_{R^*}(v)+2$ in this case), 
 one is from a doubleleton edge of $H_1^*$ and the other is a singleton edge of $H_1^*$ ($d_R(v)=d_{R^*}(v)$ in this case),    and both are singleton edges of $H_1^*$ ( $d_R(v)=d_{R^*}(v)-2$ in this case). 
 Thus $d_R(v)$ is even. 
 Since $C$ does not contain any edge of $R^*$ that is incident with $w$
 and $d_{R^*}(w)$ is even, then we know that $d_R(w)$ is also even. 
 Thus $R$ is even.

 Then  we find an Euler tour of  $R$ that begins and ends at a vertex that is not of maximum degree in the graph $H_1$ ($H_1$ is not regular as any vertex-deleted subgraph of it has at most $2k(n-2)/2-6$  edges).
 Placing the edges of  $R$ from this tour alternately in $H_{11}$ and $H_{12}$ gives the desired decomposition, each graph having maximum degree $k-1$.
 
We show that neither $H_{11}$ nor $H_{12}$ has a $(k-1)$-overfull subgraph.  Recall from before that  $H_{1}$ has at most $2k(n-2)/2 - 6$ edges in any vertex-deleted subgraph.  Therefore  $H_{1} - E(C)$ has at most $ (2k-2)(n-2)/2 - 6$ edges in any vertex-deleted subgraph. By construction, the graphs  $H_{11}$ and $H_{12}$ differ in size by at most 1, and in degree at any vertex by at most 2. Thus, arguing as before, the number of edges in any vertex-deleted subgraph of  $H_{11}$ or $H_{12}$  is at most  $(k-1)(n-2)/2-6 + 3 < (k-1)(n-2)/2 $, so that  $H_{11}$ and $H_{12}$  contain no $(k-1)$-overfull subgraphs.  It follows Theorem~\ref{thm:plantholt-shan} that  $ \chi'(H_{11}) =   \chi'(H_{12}) = k-1$.

 We now let $G_1=G-E(C)-E(H_{11})$. As  $\mu(G^{**})=2r$, we know that  each $H_i$ contains two edges between any two vertices $u,v\in V(G^{**})$
 with $e_{G^{**}}(u,v)=2r$. 
 As $H_{12}$ is a simple graph, we know that $C\cup H_{11}$ contains at least one edge  between any two vertices $u,v\in V(G^*)$
 with $e_{G^*}(u,v)=2r$.  Thus $C\cup H_{11}$ contains at least one edge  between any two vertices $u,v\in V(G)$
 with $e_{G}(u,v)=r$.  
 Thus $\mu(G_1) \le r-1$.
 
  For  $v$ of $G$,  suppose $d_{G^{*}}(v)=r k_v$ for some positive integer $k_v$. If $d_{G^*}(v)=\Delta=rk$, 
 then $d_{G_1}(v) = (r-1)k-1$. Otherwise,
 we have  $d_{G^*}(v) \le  \Delta-r$ as all the degrees  of vertices from $V(G)$  in $G^*$ are multiples of $r$. 
Then  $d_{G_1}(v) \le \max\{k_v (r-1)+1, k_v(r-1)\} $,  since the degree of $v$ in $H_{11}$ is at most $k_v+2$ and so 
 in $G_1$, the degree of $v$ is at most $rk_v -(k_v+1-2)=k_v (r-1)+1<(r-1)k-1$, as $r\ge 3$.  Thus $\Delta(G_1)=(r-1)k-1$.  
 
Finally, for any $v\in V(G_1)$, we show that $G_1-v$
contains no $\Delta(G_1)$-overfull subgraph. 
Note that $e(G_1-v) =\frac{1}{2}\sum_{i=1}^re(H_i-v)$. 
Recall from~\eqref{eqn:Hi-v} that  $e(H_i-v) \le 2k(n-2)/2-6$. Also 
$ | e(H_i-v)-e(H_j-v)|  \le 6$ for any $i,j\in [1,r]$. 
Thus $e(H_1-v) \ge \frac{1}{r} \sum_{i=1}^re(H_i-v)- 6$. 
As $e(H_{11}-v) \ge \frac{1}{2}(e(H_1-v)-(n-2))-1$,
we know that $e(H_{11}-v) \ge \frac{1}{2}(\frac{1}{r} \sum_{i=1}^re(H_i-v)-6-(n-2))-2$. Thus 
\begin{eqnarray*}
	e(G_1 - v) &=& \frac{1}{2}\sum_{i=1}^re(H_i-v)-(n-2)-e(H_{11}-v) \\
		& \le  & (\frac{1}{2}-\frac{1}{2r})\sum_{i=1}^re(H_i-v)+3+\frac{n-2}{2}+2-(n-2) \\
		&\le & (\frac{1}{2}-\frac{1}{2r}) r(2k(n-2)/2-6)-\frac{n-2}{2}+5 \\
		&=&(k(r-1)-1)(n-2)/2-3r+8 <\Delta(G_1)(n-2)/2.  
\end{eqnarray*}
Thus $G_1$ contains no $\Delta(G_1)$-overfull subgraph. By Case 1, we know that $\chi'(G_1)=\Delta(G_1)$. 
Therefore $\chi'(G) \le \chi'(G_1)+\chi'(H_{11})+\chi'(C)=\Delta(G)$
and so $\chi'(G)=\Delta(G)$. 
 \qed

\section*{Acknowledgment}
Songling Shan was supported by NSF grant DMS-2153938. 



\end{document}